\documentclass[a4paper]{article}
\pdfoutput=1
\usepackage{arxiv}

\usepackage[left=3cm,right=3cm,a4paper]{geometry}

\usepackage[utf8]{inputenc}
\usepackage[T1]{fontenc}

\usepackage[frozencache,cachedir=.]{minted}
\usepackage{zi4}
\usepackage[final,activate=true,kerning=true,protrusion=true]{microtype}
\usepackage[all]{nowidow}
\usepackage[english]{babel}
\usepackage{enumitem,dsfont,stackengine,varioref}
\usepackage[pdfpagelabels=true]{hyperref}
\usepackage{mathtools,amsfonts,amssymb}
\usepackage{cleveref,csquotes}
\usepackage[super]{nth}
\usepackage{authblk}
\usepackage{caption}
\DeclarePairedDelimiter\abs{\lvert}{\rvert}%
\makeatletter
\let\oldabs\abs
\def\abs{\@ifstar{\oldabs}{\oldabs*}}

\newtheorem{definition}{Definition}
\newtheorem{example}{Example}
\newtheorem{remark}{Remark}

\usepackage[svgnames]{xcolor}
\hypersetup{
	linktoc = page,
	pdfpagemode = UseNone,
	colorlinks,
	linkcolor={red!50!black},
	citecolor={blue!50!black},
	urlcolor={blue!80!black}
}

\usepackage{mdframed}
\usepackage{booktabs}

\newcommand{\N}{{\mathbb{N}^*}}
\newcommand{\C}{\mathbb{C}}

\newcommand{\Prime}{\mathbb{P}}
\newcommand{\one}{\mathds{1}}
\newcommand{\chapquote}[1]{\begin{quotation} \emph{#1} \end{quotation} }

\begin{document}

\title{Automated Discovery of New \texorpdfstring{$L$}{L}-Function Relations}

\author[1,*]{Hadrien Barral}
\author[1,2]{Rémi Géraud-Stewart}
\author[1]{Arthur Léonard}
\author[1]{David Naccache}
\author[1]{Quentin Vermande}
\author[1]{Samuel Vivien}
\affil[1]{Département d'informatique de l'ÉNS, \'Ecole normale supérieure, CNRS, PSL Research University, 45 rue d'Ulm, Paris, France}
\affil[2]{QPSI, Qualcomm Inc., San Diego CA, USA}
\affil[*]{Contact email: hadrien.barral@ens.fr}

\maketitle

\begin{abstract}
$L$-functions typically encode interesting information about mathematical objects. This paper reports 29 identities between such functions that hitherto never appeared in the literature. Of these we have a complete proof for 9; all others are extensively numerically checked and we welcome proofs of their (in)validity.

The method we devised to obtain these identities is a two-step process whereby a list of candidate identities is automatically generated, obtained, tested, and ultimately formally proven. The approach is however only \emph{semi-}automated as human intervention is necessary for the post-processing phase, to determine the most general form of a conjectured identity and to provide a proof for them.

This work complements other instances in the literature where automated symbolic computation has served as a
productive step toward theorem proving and can be extended in several directions further to explore the algebraic landscape of $L$-functions and similar constructions.
\end{abstract}

\keywords{$L$-functions \and conjectures \and automated}

\section*{Introduction}
Dirichlet famously introduced $L$-functions, which amongst other tools proved instrumental in establishing results in the distribution of prime numbers in infinite sequences \cite{dirichlet1889werke}. $L$-functions and their countless generalizations can be constructed for many objects, including characters\footnote{Dirichlet's original motivation}, modular forms, or elliptic curves where they are notably used to formulate the celebrated Birch--Swinnerton--Dyer conjecture \cite{hardy2013general}.

This paper focuses on $L$-functions constructed from multiplicative functions (\Cref{sec:multiplicative-functions}). The Dirichlet sums of such functions feature a particularly nice property: they can be expressed as an infinite product over the primes, such as Euler's product \cite{euler1737variae} for Riemann's $\zeta$ function (\Cref{sec:euler-product}), and are accordingly called \emph{the $L$-function's Euler product}. At the same time, the Dirichlet sum can yield a known function, such as Riemann's $\zeta$ or $\eta$ functions.
This raises the following question:

\chapquote{``Can we find \emph{remarkable} relationships between special functions (e.g., $\zeta$, logarithms etc.), or at the very least between Dirichlet sums, through the study of their Euler product?''}
Our approach consists in adapting algebraic sieving algorithms, initially designed to factor composite integers or compute discrete logarithms, to reduce the question of detecting new theorems to the finding of ``smoothness'' relationships followed by a linear algebraic processing which can be fully automated.

This method is heuristic, but the candidate identities can be tested automatically, and if they succeed, spend some time formally proving them. In doing so, we found many relations relating special functions, or at the very least Dirichlet sums. The simplest of such relations are already well-known --- see e.g., \cite{gould2000catalog} --- but we find several new, non-obvious results, which may prove useful in the studying $L$-functions. Nontrivial examples found by our algorithm are identities such as:

\begin{equation*}
\sum_{n=1}^\infty \frac{\lambda(n) \tau(n) \sigma'_2(n)}{n^6}
=\frac{\zeta(4)^2\zeta(10)\zeta(12)^2}{\zeta(6)^2\zeta(20)} = \frac{154226363 \pi^{10}}{12741871041900},
\end{equation*}
where the functions $\lambda,\tau,\sigma'$ are given hereafter.

\section{Preliminaries}

\paragraph{Notations.} We denote by $\Prime$ the set of all prime numbers and $\N$ is the set of natural numbers without $0$.

\subsection{Multiplicative functions}\label{sec:multiplicative-functions}
\begin{definition}[Multiplicative Function]
    A function $f:  \N \rightarrow \C$ is \emph{multiplicative} if for any coprime integers $x, y$, $f(xy) = f(x)f(y)$.
    We denote by $\mathcal M$ the set of multiplicative functions.
\end{definition}
\begin{example}

The functions given in Table \ref{tab:mulf} are well known to be multiplicative and are used throughout this paper. Additional multiplicative functions can be found in \cite{gould2000catalog}.

\begin{table}[!htp]
\centering
\begin{tabular}{|ccl|}\hline
$\one$      &:& $  n \mapsto 1$\\\hline
$\epsilon$&:& $ n \mapsto \begin{cases} 1 & \text{if $n = 1$} \\
            0 &  \text{otherwise}
        \end{cases}$, \emph{Kronecker $\delta_{1,n}$}\\\hline
$\text{Id}$&:&$n \mapsto n$, \emph{Identity function}\\\hline
$\varphi$&:&$n \mapsto \#\{1 \leq i \leq n:  i \wedge n = 1\}$, \emph{Euler's totient function}\\\hline
$\sigma_k$&:&$n \mapsto \sum_{i \;|\; n} i^k$, \emph{The $k^{\mbox{\scriptsize{th}}}$ divisor function}\\\hline
$\tau$&:&$\sigma_0$, the number of divisors\\\hline
$\tau_k$&:&$n \mapsto \#\{(i_1,\dotsc,i_k) \in \N^k:  \prod_{\ell=1}^k i_{\ell}  = n\}$\\
&&\emph{The number of ways to express $n$ as a product of $k$ positive factors~~~~}\\
&&\emph{Note that $\tau = \tau_2$}\\\hline
$\mu$&:& $  n \mapsto
        \begin{cases}
        (-1)^s &  \text{if } n = \prod_{i=1}^s p_i \text{ with distinct } p_1, \dotsc, p_s \in \mathbb{P} \\
        0 &  \text{otherwise}
    \end{cases}$\\
&&\emph{Möbius' function}\\\hline
$\mu_k$&:&$n \mapsto \begin{cases}
        (-1)^s &  \text{if } n = \prod_{i=1}^s p_i^k \text{ with distinct } p_1, \dotsc, p_s \in \mathbb{P} \\
        0 &  \text{otherwise}
    \end{cases}$, thus $\mu_1 = \mu$\\
&&\emph{One of the possible generalizations of Möbius' function}\\\hline

$J_k$&:&$n \mapsto \#\{(a_1, \dotsc, a_k) \in \N^k:  a_i \leq n \text{ and } (a_1, \dotsc, a_k, n) \text{ are coprime}\}$\\
&&\emph{Jordan's totient function (we have $J_k(n) = \mu(n) \ast n^k$)}\\\hline
$\lambda$&:&$  n \mapsto (-1)^r$, where $ r = \#\{(p, k) \in \mathbb{P} \times \N:  p^k \vert n\}$\\ &&\emph{Liouville's function}\\\hline
$\zeta_k$&:&$n \mapsto n^k$, where $k$ is non-negative.\\\hline
$\nu_k$&:& $ n \mapsto \begin{cases} 1 & \text{if $n$ is a $k^{\mbox{\scriptsize{th}}}$ power}\\
            0 & \text{otherwise}
        \end{cases}$\\\hline
$\xi_k$&:& $ n \mapsto \begin{cases} 1 & \text{if $n$ is $k$-free} \\
            0 &  \text{otherwise}
        \end{cases}$, \emph{where \enquote{$n$ is $k$-free} means $\forall p \in \mathbb{P}, p^k \nmid n$}\\\hline
$\theta$&:&  $n \mapsto \#\{(a,b) \in \N^2:  ab = n \text{ and } \gcd(a,b) = 1\}$\\\hline
$\sigma'_k$&:&  $n \mapsto \sum_{d \mid n} \lambda(d)d^k$, where $k$ is non-negative
\phantom{$\frac{a}{\text{{\large b}}}$} 
\\\hline
$\psi_k$ &:&$n \mapsto \sum_{d \mid n} d^k \left\lvert \mu\left(\frac{n}{d}\right) \right\rvert $\\
&&\emph{where $\psi_1$ is known as Dedekind's function} \\\hline
\end{tabular}
    \caption{Examples of multiplicative functions.}
    \label{tab:mulf}
\end{table}

\end{example}

\subsection{Dirichlet \texorpdfstring{$L$}{L}-functions}\label{sec:l-functions}
$L$-functions were formally defined and given this name by Dirichlet \cite[pp.~313--342]{dirichlet1889werke}, whose original aim was to prove that there are infinitely many primes in any (primitive) arithmetic progression.

\begin{definition}[Dirichlet $L$-functions]
If $f\in\mathcal{M}$, we define the corresponding formal series called the \emph{$L$-function} associated with $f$:
\begin{equation}
L(f, s) = \sum_{n=1}^\infty \frac{f(n)}{n^s}.
\end{equation}
\end{definition}
The well-definedness, convergence properties and analytical continuation of such sums have been extensively studied. For our purposes it is sufficient to say that if $f$ doesn't grow too fast, the corresponding $L$-functions is convergent as soon as the real part of $s$ is large enough. In particular,
\begin{align*}
L(\one, s)
& = \sum_{n=1}^\infty \frac{\one(n)}{n^s} = \sum_{n=1}^\infty \frac{1}{n^s} = \zeta(s),
\end{align*}
where $\zeta(s)$ is Riemann's zeta function \cite{riemann1859ueber}.

\begin{remark}[About convergence]
In the rest of this paper, we do not discuss the convergence of $L$-functions in detail, and assume that the formal manipulations are valid throughout. The result is that some terms in the relations we obtain may be divergent. This is only of consequence if such divergent terms end up in the final identities, at which point they are easily spotted. It is a heuristic's nature that it sometimes produces correct outputs through a reasoning that is not valid throughout. This is why we insist that the identities found by our methods, even if numerically credible, often need an independent proof to become theorems.
\end{remark}

\subsection{Euler products and Bell series}\label{sec:euler-product}
    Let $f\in\mathcal M$. Under classical convergence hypotheses, we can write an $L$-function as its Euler product:
    \begin{equation*}
    L(f, s) = \prod_{p \in \Prime} \left( \sum_{k = 0}^\infty \frac{f(p^k)}{p^{ks}} \right) \text{ written } \prod_{p \in \Prime} R_p(f, s).
    \end{equation*}
The quantity $R_p(f, s)$ is called the \emph{Bell series} associated with $f$ at $p$ and $s$ \cite[p.\,42--45]{apostol1998introduction}.

\begin{definition}[Good functions and $R$-fractions]
    Suppose $\exists R(f, s) \in \C(X)$ (the set of complex rational fractions) such that $R_p(f, s) = R(f, s)(p)$, we say that $f$ is \emph{good}, and call $R$ its $R$-fraction.
\end{definition}

\begin{example} $\varphi$ is good since:
\begin{align*}
R_p(\varphi, s) &= \sum_{k = 0}^\infty \frac{\varphi(p^k)}{p^{ks}} = 1 + \sum_{k = 1}^\infty \frac{p^k - p^{k - 1}}{p^{ks}} \\
&= 1 + \frac{p - 1}{p^s} \sum_{k = 0}^\infty \frac{1}{p^{k(s - 1)}} = 1 + \frac{p - 1}{p^s - p} \\
&= \frac{p^s - 1}{p^s - p}.
\end{align*}
\end{example}
\begin{example}
$\one$ is good since:
\begin{equation*}
R_p(\one, s) = \sum_{k = 0}^\infty \frac{1}{p^{ks}} = \frac{p^s}{p^s - 1}.
\end{equation*}
\end{example}
The key observation is that for such functions, any multiplicative relation \emph{between $R$-fractions} gives a relation \emph{between the corresponding $L$-functions}, which was one of the motivations for Bell's introduction of the eponymous series in the 1930s.
\begin{example} For $s$ large enough to guarantee convergence:
\begin{equation*}
    \frac{R(\one, s - 1)}{R(\one, s)} = \frac{X^{s - 1}}{X^{s - 1} - 1} \frac{X^s - 1}{X^s} = \frac{X^s - 1}{X^s - X} = R(\varphi, s).
\end{equation*}
Hence, we recover the well-known relation:
\begin{equation*}
    L(\varphi, s) = \frac{L(\one, s - 1)}{L(\one, s)} = \frac{\zeta(s - 1)}{\zeta(s)}.
\end{equation*}
\end{example}

\section{Adapting algebraic sieving to \texorpdfstring{$L$}{L}-functions}

In light of previous observations, we may seek multiplicative relations between $R$-fractions. We do this by adapting the algebraic sieving technique to the context of rational fractions.

\subsection{Pseudo-linear functions}
As a preliminary step, we need to obtain the $R$-fractions associated to (good) multiplicative functions. It turns out that many interesting multiplicative functions belong to a subset which is more amenable to algorithmic treatment:
\begin{definition}[Pseudo-linear function]
    $f\in \mathcal M$ is \emph{pseudo-linear} if  $\exists n \in \N$, a column vector $u \in \C(X)^n$ and a matrix $A \in M_n(\C(X))$ such that for every prime $p$,
    \begin{equation*}
        f(p^k) = \pi(A^k u)(p),
    \end{equation*}
    where $\pi$ is the projection on the first component.
We will represent such a pseudo-linear function as a pair $(A, u)$ and denote by $\mathcal P$ the set of pseudo-linear functions.
\end{definition}
\begin{example}

    For the function $\tau$ (number of divisors), we have the matrix
    $\begin{bmatrix}
    1 & 1 \\
    0 & 1 \\
    \end{bmatrix}$
    and the vector
    $\begin{bmatrix}
    1 \\
    1
    \end{bmatrix}$. We can thus see that $\tau(p^k) = k+1$.
\end{example}

\begin{remark}
    Since $R(f, s) = R(f\times \operatorname{Id}^{-s}, 0)$, we are only interested in computing $R$ when $s = 0$ for any multiplicative function $f$ such that $R(f, 0)$ is convergent.
\end{remark}

Let $f\in \mathcal P$ be represented as $(A, u)$. We have, if the sum $\sum_{k = 0}^\infty A^k(p)$ converges for every prime $p$:
\begin{align*}
    R(f, 0)(p) &= \sum_{k = 0}^\infty \pi(A^k u)(p)
    = \pi\left(\left(\sum_{k = 0}^\infty A^k(p) \right)u(p)\right) = \pi((I_n - A)^{-1}u)(p).
\end{align*}

In fact, if $(I_n - A)(p)$ is not invertible, there exists a non-zero vector $v \in \C^n$ such that $A(p)v = v$, and we have:
\[\left ( \sum_{k = 0}^\infty A^k(p) \right ) v = \sum_{k = 0}^\infty A(p)^k v = \sum_{k = 0}^\infty v\]
which contradicts the convergence of $\sum_{k = 0}^\infty A^k(p)$.
Thus, under our hypotheses, $(I_n - A)(p)$ is always invertible, and the invertibility of $I_n - A$ follows.

\subsection{Operations between \texorpdfstring{$R$}{R}-fractions}
\subsubsection{Multiplication}
Let $f, g\in \mathcal P$, represented respectively as $(A, u)$ and $(B, v)$. We have:
\begin{align*}
    fg(p^k) &= f(p^k)g(p^k) = \pi(A^k u)\pi(B^k v) \\
    &= \pi\left((A^k \otimes B^k)(u \otimes v)\right) \\
    &= \pi((A \otimes B)^k (u \otimes v))
\end{align*}
where $\otimes$ is the tensor product.
Thus, $fg$ can be represented by $(A \otimes B, u \otimes v)$.

\subsubsection{Dirichlet convolution}
If $f, g\in \mathcal P$, we denote by:
    \begin{equation*}
    f \ast g:  n \in \N \mapsto \sum_{d \;|\; n} f(n)g\left(\frac{n}{d}\right) \in \C
    \end{equation*}
    their \emph{Dirichlet convolution}. It is known that the result $f \ast g$ is again a multiplicative function \cite[Theorem 2.14]{apostol1998introduction}.
    Let $(A, u)$ and $(B, v)$ represent $f$ and $g$, respectively; we would like to compute a representation of $f \ast g$.

We have, if $A \otimes I - I \otimes B$ is invertible:
\begin{align*}
(f \ast g)(p^k) &= \sum_{i = 0}^k f(p^i)g(p^{k - i})
= \sum_{i = 0}^k \pi(A^i u)\pi(B^{k - i}v) \\
&= \sum_{i = 0}^k \pi\left((A^i \otimes B^{k - i})(u \otimes v)\right)
= \pi\left(\sum_{i = 0}^k (A^i \otimes B^{k - i})(u \otimes v)\right) \\
&= \pi((A^{k + 1} \otimes I - I \otimes B^{k + 1})(A \otimes I - I \otimes B)^{-1}(u \otimes v)) \\
&= \pi((A \otimes I)^k u') + \pi((I \otimes B)^k v') \\
&= \pi(((A \otimes I) \odot (I \otimes B))^k(u' \odot v'))
\end{align*}
where
\begin{align*}
    u' &:= (A \otimes I)(A \otimes I - I \otimes B)^{-1}(u \otimes v) \\
    v' &:= -(I \otimes B)(A \otimes I - I \otimes B)^{-1}(u \otimes v).
\end{align*}
and for any vectors $x, y$, $x \odot y := (\pi(x) + \pi(y)) \oplus x \oplus y$,
and for any matrices $A , B$:
\[A \odot B := \begin{bmatrix}
0 & A[1] & B[1] \\
0 & A & 0 \\
0 & 0 & B \\
\end{bmatrix}
\] where $M[1]$ is the first row of the matrix $M$.
Thus, $f \ast g$ can be represented by $((A \otimes I) \odot (I \otimes B), u' \odot v')$.

\begin{remark}
    If $A \otimes I - I \otimes B$ is not invertible, our implementation will detects this; and the user will have implement a representation of this multiplicative function by hand.
    A future version of our implementation could correct this problem by using a better multiplicative function representation.
\end{remark}

With this setup, it is easy to generate the rational fractions of many pseudo-linear functions: we can compute the representation of simple known multiplicative functions, then compose these representations using the product and the Dirichlet convolution, and finally compute the rational fractions.

\subsubsection{Reduction of representations}
For performance, we would like, given a representation $(A, u)$ of a pseudo-linear function $f$, to find a representation $(B, v)$ of $f$ of smaller dimension.
This can be done easily by finding a linear relation between the rows of the matrix $(A_1, \dots, A_d, u)$, then removing one of the rows using this relation.

\subsection{Generating Relations}

We can now move on to finding multiplicative relations between $R$-fractions.

\subsubsection{Holding space basis.} We introduce the following definition:
\begin{definition}[Holding space]
A set $V$ of non-zero rational fractions is a \emph{holding space} if it is a finitely generated subgroup for the multiplication.
\end{definition}
From a set $P_1, \dotsc, P_k$ of polynomials we can construct a set $\mathcal B$ which generates the same holding space, but with the property that polynomials in $\mathcal B$ are pairwise coprime. Indeed, we can proceed using the following \texttt{Insert} algorithm on $P = P_1, \dotsc, P_k$:

\begin{table}[!htp]
\centering
\setlength{\tabcolsep}{0.3em} 
\begin{tabular}{|lllll|}\hline
 \textbf{Case}  &~~~~& \textbf{Condition} & &\textbf{Perform the operation~~~}\\\hline\hline
 Case 0  &~~~~& $P = 1$ & & \texttt{Discard it}\\ \hline
 Case 1  &~~~~& $\exists Q\in \mathcal B$ such that $Q \mid P$ &~~~~&\texttt{Insert}($P/Q$)\\\hline
 Case 2  && $\exists Q\in \mathcal B$ such that $\gcd(Q,P) \neq 1$ & &$\mathcal B\phantom{'}\leftarrow \mathcal{B} -  Q$ \\
         &&                                                        & &$Q'\leftarrow Q/\gcd(Q,P)$ \\
         &&                                                        & &$\mathcal B\phantom{'}\leftarrow \mathcal B \cup Q'$ \\
         &&                                                        & &\texttt{Insert}($\gcd(Q, P)$)\\
         &&                                                        & &\texttt{Insert}($P/\gcd(Q, P)$)\\\hline
 Case 3  && otherwise                                              & &$\mathcal B\phantom{'}\leftarrow \mathcal B \cup P$ \\\hline
\end{tabular}
    \caption{\texttt{Insert} algorithm}
    \label{tab:rel_insert}
\end{table}

We convene that $\mathcal B$ is sorted so that the polynomials appearing the most frequently in the $P_i$ appear the first in $\mathcal B$.
We use this algorithm on the polynomials appearing (as numerator or denominator) in a collection of $R$-fraction, and call the resulting set $\mathcal B$ our \emph{holding space basis}.

\subsubsection{Composition matrix.} By construction, our $R$-fractions can be written as a product or ratio of elements in $\mathcal B = (b_1, \dotsc, b_r)$; furthermore this decomposition is unique up to the order of terms. Therefore, we can associate to each $R$-fraction a vector whose $i^{\mbox{\scriptsize th}}$ coefficient is the exponent of $b_i$ in this unique decomposition. Stacking these row vectors together, all our $R$-fractions give a matrix $M$.

The multiplication of two $R$-fractions corresponds to the addition of two matrix rows. More generally a multiplicative relation between $R$-fractions corresponds to a linear combination of the rows of $M$ which gives 0 --- in other terms, we are interested in finding elements in the kernel of $M$.

\section{Implementation}
This section describes data structures and other specific choices made to implement our algorithm.
The source code and all tools used for this paper are available under the GNU General Public License (\texttt{GPL-3.0-only}) license at \url{https://github.com/CrazySumsTeam/CrazySums}.

\subsection{Generation}
\label{section:generation}

The algorithm starts by generating many $R$-fractions, from which relations will be sought.
First of all, we start by giving constraints to bound the generated $L$-functions. In other words, we give to the generator the following input:
\begin{equation*}
(f_0,\dotsc,f_n),\mbox{~}(a_0,\dotsc,a_n),\mbox{~}(b_0,\dotsc,b_n),\mbox{~as well as~}(\min{}_s, \max{}_s,\max{}_{\mbox{{\scriptsize score}}})
\end{equation*}
Afterwards, the generator computes all the $L$-functions satisfying those constraints. This means for every $(j_i)_{0 \leqslant i < n}$ such that $\forall i, a_i \leqslant j_i \leqslant b_i$, we compute $R(f, s)$ such that:
\begin{align*}
f & = \prod_{i = 0}^{n-1} f_i^{j_i},  &
s & = s(f) + k, \text{ such that } \min{}_s \leqslant k \leqslant \max{}_s &
\max{}_{\mbox{{\scriptsize score}}} & \geqslant \sum_{i = 0}^{n-1} j_i
\end{align*}
where $s(f)$ is the minimal integer such that $R(f,s)$ is defined. Note that $\max_{\mbox{{\scriptsize score}}}$ is introduced to prevent generating overly complicated $L$-functions before simple ones. Not using $\max_{\mbox{{\scriptsize score}}}$ drastically increases computation time for a minimal gain in the number of interesting relations found.

One important thing to see is that if $f = \prod_i f_i^{j_i}$ has already been computed, it much faster to compute $f \cdot f_k$ than starting back from the beginning. This pruning where we just remember the last $R$-fraction computed vastly improved the computation time needed to generate the $R$-fractions. Consequently, this part wasn't modified to enable multi-threading.

\subsection{Multi-threading}

To improve performances, some computation steps have been parallelized. Those steps are the decomposition matrix computation and the basis generation.
\begin{enumerate}
    \item The decomposition matrix computation can be parallelized quite intuitively. However, as the basis has already been found, it is thus ready-only at this point. Decomposing each polynomial can then be decomposed independently of the others, making the decomposition step an embarrassingly parallel problem. Therefore, we first create a polynomial queue, on which each worker repeatedly pops a polynomial, decomposes it and stores the result in the decomposition matrix. As the actual decomposition is by far the most time-consuming step, locking the queue with simple mutexes does not create lock-contention.
    \item The basic generation step is quite trickier. We implemented this step as an iterative parallel version of the \texttt{Insert} algorithm presented in table \ref{tab:rel_insert}.  Each worker needs to modify the basis (either append a polynomial at the end or break an existing polynomial into two smaller ones). Careful use of atomic accesses, mutexes (e.g., using shared mutex locking when possible), and mathematical arguments (e.g., when a worker overtakes another one when iterating on the existing basis) were needed to avoid excessive lock contention.
\end{enumerate}

\subsection{Polynomials}

Polynomials are represented over a finite field $\mathbb{F}_p$ rather than rationals. This introduces the possibility of an error but allows for much faster operation. Errors can be made less likely by increasing $p$; however in practice, no spurious relation has been found for $p=997$.

Should spurious relations be found, running our algorithm over a range of different values of $p$ and seeking out recurring relations would filter out most, if not all, erroneous relations.

\subsection{Matrices}

The composition matrix $M$ is stored as a matrix of rationals (pairs of integers). In principle this matrix holds $n\times m$ coefficients where $n$ is the size of the basis and $m$ the number of polynomials. However, it turns out that the matrix is sparse. Hence, rather than using a 2-D array the matrix is represented as an array of rows, where a row is a vector of pairs holding the column number and the coefficient at this position. Thus, only non-zero coefficients need to be stored.

\subsection{Post-processing the results}\label{sec:post-processing-the-results}

Once the results have been computed, human intervention is necessary to turn them into mathematical statements. A typical run outputs thousands of relations, making this task quite labor-intensive. Indeed, our implementation does not infer general symbolic relations, but specific instances of them, and we did seek to display relations in a human-friendly manner. This human intervention and how we strove to ease it is discussed in more detail in \ref{apx:finding-general-symbolic-relations}.

\newpage 

\section{Results}

\subsection{New Symbolic Relations Discovered}
We give an identifier C-\texttt{XX} to the relations found using our approach, that (to the best of our knowledge) were not known before.
Known relations, and those which are special cases of others, were removed from this list. We could also directly prove the equality for some of these results, which is indicated by a \checkmark in the table (as an example, we provide the proof for C-14 in \ref{apx:example_proof}).

\begin{center}
\newcommand{\crel}[2]{
  \ifthenelse{#2=1}{\checkmark}{\phantom{\checkmark}} #1\phantom{a}
}
\def\lc{\left\lceil}
\def\rc{\right\rceil}

\renewcommand{\arraystretch}{1.15} 

\begin{tabular}{|ccc|}\hline
\textsf{ID}     & \textbf{Relation Discovered} &   \textbf{Remarks} \\\hline\hline
\crel{C-01}{1}  & $L(\sigma_n \abs{\mu}, 2n) = \zeta(n) \cdot \zeta(3n)^{-1}$&\\\hline
\crel{C-05}{1}  & ${L\left(f^2 \mu, 2n\right)} = {L\left(f \mu, n\right)}\cdot{L\left(f \abs{\mu}, n\right)}$ & for multiplicative $f$\\\hline
\crel{C-06}{0}  & ${L\left(\theta^2 \sigma_n \abs{\mu}, 2n\right)}={L\left(\theta \abs{\mu}, n\right)}^2$&\\\hline
\crel{C-07}{0}  & $L\left(\theta^2 J_n \mu, 2n\right) = {L\left(\theta \mu, n\right)}^2$ &\\\hline
\crel{C-09}{0}  & ${L\left(\phi^{n-2} \theta^{\ell} \abs{\mu}, n\right)}={L\left(\phi^{n-2} \theta^{\ell} J_2, n+2\right)}$ & part. case of C-10\\\hline
\crel{C-10}{0}  & ${L\left(\phi^i \theta^j J_k^p \abs{\mu}, i+pk\right)}={L\left(\phi^i \theta^j J_k^{p+1}, i+(p+1)k\right)}$& generalizes C-09\\\hline
\crel{C-11}{1}  & $L(J_n \mu, 2n) = \zeta(2n) \cdot \zeta(3n) \cdot \zeta(n)^{-1} \cdot \zeta(6n)^{-1}$&\\\hline
\crel{C-13}{1}  & $L\left(J_{2i} \sigma_{i+k}, 3i+2k\right) \cdot L\left(J_i \sigma_i, 2i+k\right)^{-1} =$& $i\geq1$, $k\geq2$\\
                & $\zeta(i+2k) \cdot \zeta(2i+2k)^{-1}$&\\\hline
\crel{C-14}{1}  & ${L\left(\theta \sigma_n, 2n\right)} \cdot {L\left(J_n^2, 4n\right)}^{-1} = {\zeta\left(n\right)^2}$&\\\hline
\crel{C-15}{0}  & ${{L\left(J_n \sigma_n, 3n\right)}\cdot{L\left(J_n^2, 4n\right)}}\cdot{{L\left(J_n^2 \sigma_n^2, 5n\right)}^{-1}}=$& \\
                & ${{\zeta\left(2n\right)}}\cdot{{\zeta\left(3n\right)}^{-1}}$&\\\hline
\crel{C-16}{0}  & $J_{2n} \mu$ = $J_n \sigma_n \mu$& \checkmark\phantom{a} for $n=1$\\
                & (Note: \checkmark\phantom{a} $\forall n, J_{2n} \neq J_n \sigma_n$)&\\\hline
\crel{C-18}{0}  & $L\left(\tau \sigma_k, n\right) = \zeta(n-k)^2 \cdot \zeta(n)^2 \cdot \zeta(2n-k)^{-1}$&\\\hline
\crel{C-19}{0}  & ${{L\left(\lambda \tau \sigma_k, n\right)}}={{\zeta\left(2n-2k\right)}^2\cdot{\zeta\left(2n\right)}^2}\cdot$&\\
                & ${\zeta\left(n-k\right)}^{-2}\cdot{\zeta\left(n\right)}^{-2}\cdot{\zeta\left(2n-k\right)}^{-1}$&\\\hline
\crel{C-20}{0}  & $L\left(\tau J_n, 2n\right) \cdot L\left(J_n^2, 3n\right)^{-1} = \zeta(n)$&$n\geq3$\\\hline
\crel{C-21}{0}  & ${L\left(\tau \theta, n\right)} \cdot {L\left(\tau J_n^2, 3n\right)}^{-1} = {\zeta(n)^2}$&$n\geq3$\\\hline
\crel{C-22}{0}  & ${L\left(\lambda \tau_k, n\right)}={{\zeta(2n)}^k} \cdot {\zeta(n)}^{-k}$&\checkmark\phantom{a} for $k=2$\\\hline
\crel{C-23a}{1} & $\lambda \nu_{2k} = 1$&trivial\\
\crel{C-23b}{1} & $\lambda \nu_{2k+1} = \lambda$&trivial\\\hline
\crel{C-24}{0}  & ${{L\left(\tau \sigma'_k, n\right)}}={{\zeta\left(n\right)}^2\cdot{\zeta\left(2n-2k\right)}^2\cdot{\zeta\left(2n-k\right)}}\cdot$&\\
                & ${\zeta\left(n-k\right)}^{-2}\cdot{\zeta\left(4n-2k\right)}^{-1}$&\\\hline
\crel{C-25}{0}  & ${{L\left(\lambda \tau \sigma'_k, n\right)}}={{\zeta\left(2n\right)}^2\cdot{\zeta\left(n-k\right)}^2\cdot{\zeta\left(2n-k\right)}}\cdot$&\\
                & ${\zeta\left(n\right)}^{-2}\cdot{\zeta\left(4n-2k\right)}^{-1}$&\\\hline
\crel{C-26}{1}  & $\xi_k\xi_{\ell} = \xi_{\min(k,\ell)}$&trivial\\\hline
\crel{C-27}{0}  & ${L\left(\lambda \xi_k, n\right)}={\zeta(2n)}\cdot{\zeta(n)^{-1}}\cdot{\zeta(kn)^{-1}}$&$k$ odd\\\hline
\crel{C-28}{0}  & ${L\left(\lambda \xi_k, n\right)} = {\zeta(2n)}\cdot{\zeta(kn)}\cdot{\zeta(n)}^{-1}\cdot{\zeta(2kn)}^{-1}$&$k$ even\\\hline
\crel{C-29}{0}  & ${L\left(m \mapsto J_k(m^\ell), n\ell\right)}={\zeta\left((n-k)\ell\right)} \cdot {\zeta\left((n-k)\ell+k\right)}^{-1}$&\\\hline
\crel{C-30}{0}  & $L\left(m \mapsto \theta(m^\ell ), n\ell\right) = {{\zeta\left(\ell n\right)}^2} \cdot {\zeta\left(2\ell n\right)}^{-1}$&\\\hline
\crel{C-31}{1}  & $L\left(m \mapsto \xi_k(m^\ell), n\ell\right) = 1$&$k\leq\ell$, trivial\\\hline
\crel{C-32}{0}  & $L\left(m \mapsto \xi_k(m^\ell), n\ell\right) = \zeta(\ell n) \cdot \zeta\left(\lc\frac{k}{\ell}\rc \ell n\right)^{-1} $&$k>\ell$\\\hline
\crel{C-33}{0}  & $L\left(\tau \xi_2, n\right) = L\left(\theta \xi_2, n\right)$&\\\hline
\crel{C-34}{0}  & $L\left(m \mapsto \left(\theta J_{k}\right)(m^\ell), n\ell\right) = L\left(\theta J_{k}, n\ell-k(\ell-1)\right)$&\\\hline
\end{tabular}
\end{center}

\subsection{Already-Known results}
Unsurprisingly, amongst the many identities found through our algorithm, some were already known. In the relation catalog \cite{gould2000catalog}, our algorithm has been able to find the following known symbolic relations:
\begin{center}
    D-[2-6,9-13,15,18,21,22,24-28,30,37-43,46,47,49-53,55,58]
\end{center}
\begin{remark}
Relations D-27 and D-51 were found to be erroneous in \cite{gould2000catalog}, but the original source \cite{mcccarthy} used by \cite{gould2000catalog} is correct and consistent with our results.
\end{remark}

\begin{remark}
Although our algorithm automatically identifies relations, its output is not in mathematical form and human intervention was necessary to infer the general forms of the relations and formulate the discovered conjectures under a parameterized form. As a result,
a (large) part of the program's output has not yet been processed. This task could maybe be automated, at least partly, and several other intriguing results still await to be uncovered.
\end{remark}

\begin{remark}
Obtaining even more relations would be possible, by investing more computational power and/or modifying the implementation to explore different branches, with the same caveat as above that turning the output into mathematical results is a labor-intensive task.
\end{remark}

\subsection{Performance}

It is unknown how many relations exist in the search space. We report on experimental results based on repeated runs, with different  configurations to find different types of relations (see section \ref{section:generation} about generation). All runs were measured on an 8 cores Intel i7-10700 CPU with 16GB of RAM.
\begin{itemize}
\item The default \emph{small} configuration of our tool runs in 1.2 seconds. It generates 1858 $R$-fractions, leading to a basis of 1517 polynomials. 659 relations are found.
\item An extended, \emph{large} configuration setup runs in 17 minutes. It generates 55198 $R$-fractions, leading to a basis of 52667 polynomials. 5590 relations are found.
\end{itemize}
We provide in \Cref{tab:run_default_breakdown} the breakdown of the relations found by those runs.

\begin{table}[!htp]
    \centering
    \begin{tabular}{|ll|c|c|}\hline
     \textbf{Classification} &~~~& \textbf{Small configuration} & \textbf{Large configuration}\\\hline\hline
     Known relations         & & $513$ & $1305$ \\ \hline
     Unclassified relations  & &  $19$ & $3923$ \\ \hline
     New relations           & & $127$ &  $362$ \\ \hline
     \textbf{Total}          & & $\mathbf{659}$ & $\mathbf{5590}$ \\ \hline
    \end{tabular}
    \caption{Breakdown of the number of relations found by our tool in different configurations}
    \label{tab:run_default_breakdown}
\end{table}

\begin{remark}
There is no precise sense in which the running time of these algorithms can be analyzed, in part because such an analysis would require having \enquote{density} results on the distribution of $L$-functions (in the usual cryptanalytic context of these methods, it is e.g., the distribution of prime numbers that matters, which is of course much better understood).
\end{remark}

\subsection{Future work: Towards new functions}
The techniques described in this article call for deeper investigations and open new research horizons. We give a couple of research directions that we think could extend the breadth of automated theorem discovery.

The sieving framework can be used beyond our use in this paper, to identify further properties of $L$-functions, and possibly other functions. Naturally, the holding space can be extended with more functions. Still, there is a limit to how many interesting multiplicative functions we can find, and it makes sense to instead extend the kind of relations that we look for.

Indeed, we may let go of the desire that $R_p(f, s)$ be expressed as a polynomial, extending the realm of relations we find. We can also consider $L$-functions for which the Bell series are not rational fractions: if it is a product of a rational fraction and some other function, say a logarithm of a rational fraction, we could use algebraic sieving on those other functions.

A further visual optimization would consist in adding to our algorithm known ``shortcut relationships'' to eliminate before the linear algebra step terms such as:
\begin{equation*}
    \prod_{p\, \in\, \mathbb{P}} \Big(\frac{p^{2s}+1}{p^{2s}-1}\Big) = \frac{B_{2s}^2 \times (4s)!}{2
B_{4s} (2s)!^2}.
\end{equation*}
This would avoid several $\zeta(x)$ in the resulting expressions but will not fundamentally modify the discovered relations. Here $B_n$ denotes the absolute value of the $n^{\mbox{\scriptsize th}}$ Bernoulli number.

We also limited ourselves to the convergent setting, but there are conceivably some identities that involve non-convergent functions such as:
\begin{equation*}
\begin{aligned}
&\lim_{m\rightarrow \infty}
\frac{m^7}{\log m}\cdot
\frac{\sum_{i=1}^m \sigma(i)^4 \phi(i)^2i^{-3}}{ \sum_{i=1}^m \sigma(i)\phi(i)}=1
\mbox{~~or~~} \\
&\lim_{m\rightarrow \infty}\frac{m \sum_{i=1}^m \sigma(i)^2 \phi(i)^2 i^{-3}}{(\sum_{i=1}^m \sigma(i)\phi(i))(\sum_{i=1}^m\phi(i)^2 i^{-4})}=\frac{\zeta(3)}{\zeta{(2)}}.
\end{aligned}
\end{equation*}
How to derive such identities in a mathematically consistent (and automated?) way?

\section*{Acknowledgments}

We would like to thank Éric Brier for proofreading an early version of this work and providing useful feedback.

\bibliographystyle{alpha}
\bibliography{crazy}

\newpage 

\appendix
\section{Example proof of a relation}
\label{apx:example_proof}

We have not proved all the relations, but we checked some of them by hand. This process could be automated in the future.
Let us show how to prove a relation.

For example, we would like to prove (C-14) that:
\begin{equation*}
    \frac{{L\left(\theta \sigma_n, 2n\right)}}{{L\left(J_n^2, 4n\right)}}={{\zeta^2\left(n\right)}}
\end{equation*}
We already proved that:
\begin{equation*}
R(\one, n) = \frac{X^n}{X^n - 1}
\end{equation*}
We have:
\begin{align*}
R_p(\theta\sigma_n, 2n) &= \sum_{k = 0}^\infty \frac{\theta(p^k)\sigma_n(p^k)}{p^{2nk}}
= 1 + \sum_{k = 1}^\infty \frac{2 \sum_{i = 0}^k p^{in}}{p^{2nk}}
= 1 + 2 \sum_{k = 1}^\infty \frac{p^{(k + 1)n} - 1}{(p^n - 1)p^{2nk}} \\
&= 1 + \frac2{p^n - 1} \sum_{k = 1}^\infty \frac{p^{(k + 1)n} - 1}{p^{2nk}}
= 1 + \frac2{p^n - 1} \sum_{k = 1}^\infty \frac{1}{p^{(k - 1)n}} - \frac{1}{p^{2nk}} \\
&= 1 + \frac{2p^n}{(p^n - 1)^2} - \frac2{(p^n - 1)(p^{2n} - 1)} \\
&= \frac{p^{4n} - 2p^n + 1}{(p^n - 1)^2(p^{2n} - 1)}
\end{align*}
Thus:
\begin{equation*}
R(\theta\sigma_n, 2n) = \frac{X^{4n} - 2X^n + 1}{(X^n - 1)^2(X^{2n} - 1)}
\end{equation*}
Also:
\begin{align*}
R_p(J_n^2, 4n) &= \sum_{k = 0}^\infty \frac{J_n^2(p^k)}{p^{4nk}} = 1 + \sum_{k = 1}^\infty \frac{(p^{nk} - p^{n(k-1)})^2}{p^{4nk}} \\
&= 1 + (1 - 2p^{-n} + p^{-2n})\sum_{k = 1}^\infty \frac{1}{p^{2nk}} \\
&= 1 + \frac{p^{2n} - 2p^{n} + 1}{p^{4n} - p^{2n}} \\
&= \frac{p^{4n} - 2p^n + 1}{p^{4n} - p^{2n}}
\end{align*}
Thus:
\begin{equation*}
    R(J_n^2, 4n) = \frac{X^{4n} - 2X^n + 1}{X^{4n} - X^{2n}}
\end{equation*}
Now:
\begin{align*}
    \frac{R(\theta\sigma_n, 2n)}{R(J_n^2, 4n)} &= \frac{X^{4n} - 2X^n + 1}{(X^n - 1)^2(X^{2n} - 1)}\frac{X^{4n} - X^{2n}}{X^{4n} - 2X^n + 1}
    = \frac{X^{2n}}{(X^n - 1)^2} \\
    & = R(\one, n)^2
\end{align*}
Thus, under convergence hypotheses, we have the expected relation.

\section{From relations instances to general symbolic formulae}
\label{apx:finding-general-symbolic-relations}

As mentioned in subsection \ref{sec:post-processing-the-results}, human intervention is necessary to turn the output of our tool into general symbolic mathematical statements.

Significant software engineering work was needed to automatically sort, classify and display relations to simplify human intervention. As a rough approximation, half of the source code is dedicated to this task.

\subsection{Pretty-printer}
We hence implemented a pretty-printer which generates two different output formats (both methods are implemented together to prevent code duplication):
\begin{itemize}
    \item The first format uses Unicode characters for mathematical symbols (such as $\sigma$) and dumps the formulae in the shell. This enables the user to directly see the result once computed and analyze them on the flight.
    \item The second format is \LaTeX{} code.
\end{itemize}

\subsection{Formula classifier}
Moreover, a formula classifier has been implemented to isolate known finds from unknown ones.
For example, ${{L\left(\lambda \sigma'_k, s\right)}}={{\zeta\left(2s\right)}\cdot{\zeta\left(s-k\right)}\cdot{\zeta\left(s\right)}^{-1}}$ (D-25) is classified with the following definition:
\begin{minted}[fontsize=\footnotesize]{c++}
    vector<pair<HFormula, Rational>> d_25 {
      {HFormulaLFunction(HFormulaProduct(
        HFormulaLeaf(FormulaNode::LEAF_LIOUVILLE),
        HFormulaLeaf(FormulaNode::LEAF_SIGMA_PRIME,
          (FormulaNode::LeafExtraArg){.k = FormulaNode::Symbolic("k"), .l = 0})
       ), FormulaNode::Symbolic("s")), Rational(1)},
      {HFormulaLFunction(HFormulaOne(), FormulaNode::Symbolic("s")), Rational(1)},
      {HFormulaLFunction(HFormulaOne(), FormulaNode::Symbolic("s+-k")), Rational(-1)},
      {HFormulaLFunction(HFormulaOne(), FormulaNode::Symbolic("2*s")), Rational(-1)},
    };
\end{minted}

\subsection{Practical example}
As an example\footnote{This example is artificial, but it was necessary to reduce the size of the tool output to ease the reading comprehension. Typical runs have $100$ to $10000$ output lines.}, let us generate formulae of the form $L\left(\lambda \tau^i {\sigma'}_{2}^{j}, s\right), 0\leq{}i\leq1, 0\leq{}j\leq2, s=s(f)+k, 0\leq{}k\leq2$. Having populated our tools with known formulae from \cite{gould2000catalog}, the output if the following:

\noindent
\hspace*{1cm}{\small\texttt{{\color{DarkGreen}[D-25]} L{\color{Purple}(}$\lambda$ $\sigma$', 3{\color{Purple})} {\color{FireBrick}=} {\color{Purple}(}$\zeta${\color{Purple}(}2{\color{Purple})} $\zeta${\color{Purple}(}6{\color{Purple})}{\color{Purple})} {\color{DarkGreen}/} $\zeta${\color{Purple}(}3{\color{Purple})}}}\\
\hspace*{1cm}{\small\texttt{{\color{DarkGreen}[D-25]} L{\color{Purple}(}$\lambda$ $\sigma$', 4{\color{Purple})} {\color{FireBrick}=} {\color{Purple}(}$\zeta${\color{Purple}(}3{\color{Purple})} $\zeta${\color{Purple}(}8{\color{Purple})}{\color{Purple})} {\color{DarkGreen}/} $\zeta${\color{Purple}(}4{\color{Purple})}}}\\
\hspace*{1cm}{\small\texttt{{\color{DarkGreen}[D-25]} L{\color{Purple}(}$\lambda$ $\sigma$', 5{\color{Purple})} {\color{FireBrick}=} {\color{Purple}(}$\zeta${\color{Purple}(}4{\color{Purple})} $\zeta${\color{Purple}(}10{\color{Purple})}{\color{Purple})} {\color{DarkGreen}/} $\zeta${\color{Purple}(}5{\color{Purple})}}}\\
\hspace*{1cm}{\small\texttt{{\color{DarkGreen}[D-42]} L{\color{Purple}(}$\lambda$ $\sigma$'\^{}2, 4{\color{Purple})} {\color{FireBrick}=} {\color{Purple}(}$\zeta${\color{Purple}(}3{\color{Purple})}\^{}2 $\zeta${\color{Purple}(}8{\color{Purple})}{\color{Purple})} {\color{DarkGreen}/} {\color{Purple}(}$\zeta${\color{Purple}(}2{\color{Purple})} $\zeta${\color{Purple}(}6{\color{Purple})}{\color{Purple})}}}\\
\hspace*{1cm}{\small\texttt{{\color{DarkGreen}[D-42]} L{\color{Purple}(}$\lambda$ $\sigma$'\^{}2, 5{\color{Purple})} {\color{FireBrick}=} {\color{Purple}(}$\zeta${\color{Purple}(}4{\color{Purple})}\^{}2 $\zeta${\color{Purple}(}6{\color{Purple})} $\zeta${\color{Purple}(}10{\color{Purple})}{\color{Purple})} {\color{DarkGreen}/} {\color{Purple}(}$\zeta${\color{Purple}(}3{\color{Purple})} $\zeta${\color{Purple}(}5{\color{Purple})} $\zeta${\color{Purple}(}8{\color{Purple})}{\color{Purple})}}}\\
\hspace*{1cm}{\small\texttt{{\color{DarkGreen}[D-42]} L{\color{Purple}(}$\lambda$ $\sigma$'\^{}2, 6{\color{Purple})} {\color{FireBrick}=} {\color{Purple}(}$\zeta${\color{Purple}(}5{\color{Purple})}\^{}2 $\zeta${\color{Purple}(}8{\color{Purple})} $\zeta${\color{Purple}(}12{\color{Purple})}{\color{Purple})} {\color{DarkGreen}/} {\color{Purple}(}$\zeta${\color{Purple}(}4{\color{Purple})} $\zeta${\color{Purple}(}6{\color{Purple})} $\zeta${\color{Purple}(}10{\color{Purple})}{\color{Purple})}}}\\
\hspace*{1cm}{\small\texttt{{\color{DarkGreen}[D-53]} L{\color{Purple}(}$\lambda$, 2{\color{Purple})} {\color{FireBrick}=} $\zeta${\color{Purple}(}4{\color{Purple})} {\color{DarkGreen}/} $\zeta${\color{Purple}(}2{\color{Purple})}}}\\
\hspace*{1cm}{\small\texttt{{\color{DarkGreen}[D-53]} L{\color{Purple}(}$\lambda$, 3{\color{Purple})} {\color{FireBrick}=} $\zeta${\color{Purple}(}6{\color{Purple})} {\color{DarkGreen}/} $\zeta${\color{Purple}(}3{\color{Purple})}}}\\
\hspace*{1cm}{\small\texttt{{\color{DarkGreen}[D-53]} L{\color{Purple}(}$\lambda$, 4{\color{Purple})} {\color{FireBrick}=} $\zeta${\color{Purple}(}8{\color{Purple})} {\color{DarkGreen}/} $\zeta${\color{Purple}(}4{\color{Purple})}}}\\
\hspace*{1cm}{\small\texttt{{\color{FireBrick}[!!!!]} L{\color{Purple}(}$\lambda$ $\tau$, 4{\color{Purple})} {\color{FireBrick}=} $\zeta${\color{Purple}(}8{\color{Purple})}\^{}2 {\color{DarkGreen}/} $\zeta${\color{Purple}(}4{\color{Purple})}\^{}2}}\\
\hspace*{1cm}{\small\texttt{{\color{FireBrick}[!!!!]} L{\color{Purple}(}$\lambda$ $\tau$, 5{\color{Purple})} {\color{FireBrick}=} $\zeta${\color{Purple}(}10{\color{Purple})}\^{}2 {\color{DarkGreen}/} $\zeta${\color{Purple}(}5{\color{Purple})}\^{}2}}\\
\hspace*{1cm}{\small\texttt{{\color{FireBrick}[!!!!]} L{\color{Purple}(}$\lambda$ $\tau$, 6{\color{Purple})} {\color{FireBrick}=} $\zeta${\color{Purple}(}12{\color{Purple})}\^{}2 {\color{DarkGreen}/} $\zeta${\color{Purple}(}6{\color{Purple})}\^{}2}}\\
\hspace*{1cm}{\small\texttt{{\color{FireBrick}[!!!!]} L{\color{Purple}(}$\lambda$ $\tau$ $\sigma$', 5{\color{Purple})} {\color{FireBrick}=} {\color{Purple}(}$\zeta${\color{Purple}(}4{\color{Purple})}\^{}2 $\zeta${\color{Purple}(}9{\color{Purple})} $\zeta${\color{Purple}(}10{\color{Purple})}\^{}2{\color{Purple})} {\color{DarkGreen}/} {\color{Purple}(}$\zeta${\color{Purple}(}5{\color{Purple})}\^{}2 $\zeta${\color{Purple}(}18{\color{Purple})}{\color{Purple})}}}\\
\hspace*{1cm}{\small\texttt{{\color{FireBrick}[!!!!]} L{\color{Purple}(}$\lambda$ $\tau$ $\sigma$', 6{\color{Purple})} {\color{FireBrick}=} {\color{Purple}(}$\zeta${\color{Purple}(}5{\color{Purple})}\^{}2 $\zeta${\color{Purple}(}11{\color{Purple})} $\zeta${\color{Purple}(}12{\color{Purple})}\^{}2{\color{Purple})} {\color{DarkGreen}/} {\color{Purple}(}$\zeta${\color{Purple}(}6{\color{Purple})}\^{}2 $\zeta${\color{Purple}(}22{\color{Purple})}{\color{Purple})}}}\\
\hspace*{1cm}{\small\texttt{{\color{FireBrick}[!!!!]} L{\color{Purple}(}$\lambda$ $\tau$ $\sigma$', 7{\color{Purple})} {\color{FireBrick}=} {\color{Purple}(}$\zeta${\color{Purple}(}6{\color{Purple})}\^{}2 $\zeta${\color{Purple}(}13{\color{Purple})} $\zeta${\color{Purple}(}14{\color{Purple})}\^{}2{\color{Purple})} {\color{DarkGreen}/} {\color{Purple}(}$\zeta${\color{Purple}(}7{\color{Purple})}\^{}2 $\zeta${\color{Purple}(}26{\color{Purple})}{\color{Purple})}}}

\noindent The human step is then to conjecture that the 3 last lines form a pattern, namely:
$${{L\left(\lambda \tau \sigma'_k, s\right)}}={{\zeta\left(2s\right)}^2\cdot{\zeta\left(s-k\right)}^2\cdot{\zeta\left(2s-k\right)}\cdot{\zeta\left(s\right)^{-2}}\cdot{\zeta\left(4s-2k\right)}^{-1}}$$.
To strengthen evidence of this conjecture, we re-run the tool with a tighter formulae generation, but with a larger $\text{max}_s$.

\noindent
\hspace*{1cm}{\small\texttt{{\color{FireBrick}[!!!!]} L{\color{Purple}(}$\lambda$ $\tau$ $\sigma$', 5{\color{Purple})} {\color{FireBrick}=} {\color{Purple}(}$\zeta${\color{Purple}(}4{\color{Purple})}\^{}2 $\zeta${\color{Purple}(}9{\color{Purple})} $\zeta${\color{Purple}(}10{\color{Purple})}\^{}2{\color{Purple})} {\color{DarkGreen}/} {\color{Purple}(}$\zeta${\color{Purple}(}5{\color{Purple})}\^{}2 $\zeta${\color{Purple}(}18{\color{Purple})}{\color{Purple})}}}\\
\hspace*{1cm}{\small\texttt{{\color{FireBrick}[!!!!]} L{\color{Purple}(}$\lambda$ $\tau$ $\sigma$', 6{\color{Purple})} {\color{FireBrick}=} {\color{Purple}(}$\zeta${\color{Purple}(}5{\color{Purple})}\^{}2 $\zeta${\color{Purple}(}11{\color{Purple})} $\zeta${\color{Purple}(}12{\color{Purple})}\^{}2{\color{Purple})} {\color{DarkGreen}/} {\color{Purple}(}$\zeta${\color{Purple}(}6{\color{Purple})}\^{}2 $\zeta${\color{Purple}(}22{\color{Purple})}{\color{Purple})}}}\\
\hspace*{1cm}{\small\texttt{{\color{FireBrick}[!!!!]} L{\color{Purple}(}$\lambda$ $\tau$ $\sigma$', 7{\color{Purple})} {\color{FireBrick}=} {\color{Purple}(}$\zeta${\color{Purple}(}6{\color{Purple})}\^{}2 $\zeta${\color{Purple}(}13{\color{Purple})} $\zeta${\color{Purple}(}14{\color{Purple})}\^{}2{\color{Purple})} {\color{DarkGreen}/} {\color{Purple}(}$\zeta${\color{Purple}(}7{\color{Purple})}\^{}2 $\zeta${\color{Purple}(}26{\color{Purple})}{\color{Purple})}}}\\
\hspace*{1cm}{\small\texttt{{\color{FireBrick}[!!!!]} L{\color{Purple}(}$\lambda$ $\tau$ $\sigma$', 8{\color{Purple})} {\color{FireBrick}=} {\color{Purple}(}$\zeta${\color{Purple}(}7{\color{Purple})}\^{}2 $\zeta${\color{Purple}(}15{\color{Purple})} $\zeta${\color{Purple}(}16{\color{Purple})}\^{}2{\color{Purple})} {\color{DarkGreen}/} {\color{Purple}(}$\zeta${\color{Purple}(}8{\color{Purple})}\^{}2 $\zeta${\color{Purple}(}30{\color{Purple})}{\color{Purple})}}}\\
\hspace*{1cm}{\small\texttt{{\color{FireBrick}[!!!!]} L{\color{Purple}(}$\lambda$ $\tau$ $\sigma$', 9{\color{Purple})} {\color{FireBrick}=} {\color{Purple}(}$\zeta${\color{Purple}(}8{\color{Purple})}\^{}2 $\zeta${\color{Purple}(}17{\color{Purple})} $\zeta${\color{Purple}(}18{\color{Purple})}\^{}2{\color{Purple})} {\color{DarkGreen}/} {\color{Purple}(}$\zeta${\color{Purple}(}9{\color{Purple})}\^{}2 $\zeta${\color{Purple}(}34{\color{Purple})}{\color{Purple})}}}\\
\hspace*{1cm}{\small\texttt{{\color{FireBrick}[!!!!]} L{\color{Purple}(}$\lambda$ $\tau$ $\sigma$', 10{\color{Purple})} {\color{FireBrick}=} {\color{Purple}(}$\zeta${\color{Purple}(}9{\color{Purple})}\^{}2 $\zeta${\color{Purple}(}19{\color{Purple})} $\zeta${\color{Purple}(}20{\color{Purple})}\^{}2{\color{Purple})} {\color{DarkGreen}/} {\color{Purple}(}$\zeta${\color{Purple}(}10{\color{Purple})}\^{}2 $\zeta${\color{Purple}(}38{\color{Purple})}{\color{Purple})}}}\\
\hspace*{1cm}{\small\texttt{{\color{FireBrick}[!!!!]} L{\color{Purple}(}$\lambda$ $\tau$ $\sigma$', 11{\color{Purple})} {\color{FireBrick}=} {\color{Purple}(}$\zeta${\color{Purple}(}10{\color{Purple})}\^{}2 $\zeta${\color{Purple}(}21{\color{Purple})} $\zeta${\color{Purple}(}22{\color{Purple})}\^{}2{\color{Purple})} {\color{DarkGreen}/} {\color{Purple}(}$\zeta${\color{Purple}(}11{\color{Purple})}\^{}2 $\zeta${\color{Purple}(}42{\color{Purple})}{\color{Purple})}}}\\
\hspace*{1cm}{\small\texttt{{\color{FireBrick}[!!!!]} L{\color{Purple}(}$\lambda$ $\tau$ $\sigma$', 12{\color{Purple})} {\color{FireBrick}=} {\color{Purple}(}$\zeta${\color{Purple}(}11{\color{Purple})}\^{}2 $\zeta${\color{Purple}(}23{\color{Purple})} $\zeta${\color{Purple}(}24{\color{Purple})}\^{}2{\color{Purple})} {\color{DarkGreen}/} {\color{Purple}(}$\zeta${\color{Purple}(}12{\color{Purple})}\^{}2 $\zeta${\color{Purple}(}46{\color{Purple})}{\color{Purple})}}}\\
\hspace*{1cm}{\small\texttt{{\color{FireBrick}[!!!!]} L{\color{Purple}(}$\lambda$ $\tau$ $\sigma$', 13{\color{Purple})} {\color{FireBrick}=} {\color{Purple}(}$\zeta${\color{Purple}(}12{\color{Purple})}\^{}2 $\zeta${\color{Purple}(}25{\color{Purple})} $\zeta${\color{Purple}(}26{\color{Purple})}\^{}2{\color{Purple})} {\color{DarkGreen}/} {\color{Purple}(}$\zeta${\color{Purple}(}13{\color{Purple})}\^{}2 $\zeta${\color{Purple}(}50{\color{Purple})}{\color{Purple})}}}

\noindent Now having sufficient evidence, we give a name to this conjecture (here C-25), and process to add it to our classifier.

The new output of our tool allows us to focus on the remaining unknown finds:

\noindent
\hspace*{1cm}{\small\texttt{{\color{DarkGreen}[D-25]} L{\color{Purple}(}$\lambda$ $\sigma$', 3{\color{Purple})} {\color{FireBrick}=} {\color{Purple}(}$\zeta${\color{Purple}(}2{\color{Purple})} $\zeta${\color{Purple}(}6{\color{Purple})}{\color{Purple})} {\color{DarkGreen}/} $\zeta${\color{Purple}(}3{\color{Purple})}}}\\
\hspace*{1cm}{\small\texttt{{\color{DarkGreen}[D-25]} L{\color{Purple}(}$\lambda$ $\sigma$', 4{\color{Purple})} {\color{FireBrick}=} {\color{Purple}(}$\zeta${\color{Purple}(}3{\color{Purple})} $\zeta${\color{Purple}(}8{\color{Purple})}{\color{Purple})} {\color{DarkGreen}/} $\zeta${\color{Purple}(}4{\color{Purple})}}}\\
\hspace*{1cm}{\small\texttt{{\color{DarkGreen}[D-25]} L{\color{Purple}(}$\lambda$ $\sigma$', 5{\color{Purple})} {\color{FireBrick}=} {\color{Purple}(}$\zeta${\color{Purple}(}4{\color{Purple})} $\zeta${\color{Purple}(}10{\color{Purple})}{\color{Purple})} {\color{DarkGreen}/} $\zeta${\color{Purple}(}5{\color{Purple})}}}\\
\hspace*{1cm}{\small\texttt{{\color{DarkGreen}[D-42]} L{\color{Purple}(}$\lambda$ $\sigma$'\^{}2, 4{\color{Purple})} {\color{FireBrick}=} {\color{Purple}(}$\zeta${\color{Purple}(}3{\color{Purple})}\^{}2 $\zeta${\color{Purple}(}8{\color{Purple})}{\color{Purple})} {\color{DarkGreen}/} {\color{Purple}(}$\zeta${\color{Purple}(}2{\color{Purple})} $\zeta${\color{Purple}(}6{\color{Purple})}{\color{Purple})}}}\\
\hspace*{1cm}{\small\texttt{{\color{DarkGreen}[D-42]} L{\color{Purple}(}$\lambda$ $\sigma$'\^{}2, 5{\color{Purple})} {\color{FireBrick}=} {\color{Purple}(}$\zeta${\color{Purple}(}4{\color{Purple})}\^{}2 $\zeta${\color{Purple}(}6{\color{Purple})} $\zeta${\color{Purple}(}10{\color{Purple})}{\color{Purple})} {\color{DarkGreen}/} {\color{Purple}(}$\zeta${\color{Purple}(}3{\color{Purple})} $\zeta${\color{Purple}(}5{\color{Purple})} $\zeta${\color{Purple}(}8{\color{Purple})}{\color{Purple})}}}\\
\hspace*{1cm}{\small\texttt{{\color{DarkGreen}[D-42]} L{\color{Purple}(}$\lambda$ $\sigma$'\^{}2, 6{\color{Purple})} {\color{FireBrick}=} {\color{Purple}(}$\zeta${\color{Purple}(}5{\color{Purple})}\^{}2 $\zeta${\color{Purple}(}8{\color{Purple})} $\zeta${\color{Purple}(}12{\color{Purple})}{\color{Purple})} {\color{DarkGreen}/} {\color{Purple}(}$\zeta${\color{Purple}(}4{\color{Purple})} $\zeta${\color{Purple}(}6{\color{Purple})} $\zeta${\color{Purple}(}10{\color{Purple})}{\color{Purple})}}}\\
\hspace*{1cm}{\small\texttt{{\color{DarkGreen}[D-53]} L{\color{Purple}(}$\lambda$, 2{\color{Purple})} {\color{FireBrick}=} $\zeta${\color{Purple}(}4{\color{Purple})} {\color{DarkGreen}/} $\zeta${\color{Purple}(}2{\color{Purple})}}}\\
\hspace*{1cm}{\small\texttt{{\color{DarkGreen}[D-53]} L{\color{Purple}(}$\lambda$, 3{\color{Purple})} {\color{FireBrick}=} $\zeta${\color{Purple}(}6{\color{Purple})} {\color{DarkGreen}/} $\zeta${\color{Purple}(}3{\color{Purple})}}}\\
\hspace*{1cm}{\small\texttt{{\color{DarkGreen}[D-53]} L{\color{Purple}(}$\lambda$, 4{\color{Purple})} {\color{FireBrick}=} $\zeta${\color{Purple}(}8{\color{Purple})} {\color{DarkGreen}/} $\zeta${\color{Purple}(}4{\color{Purple})}}}\\
\hspace*{1cm}{\small\texttt{{\color{Blue}[C-25]} L{\color{Purple}(}$\lambda$ $\tau$ $\sigma$', 5{\color{Purple})} {\color{FireBrick}=} {\color{Purple}(}$\zeta${\color{Purple}(}4{\color{Purple})}\^{}2 $\zeta${\color{Purple}(}9{\color{Purple})} $\zeta${\color{Purple}(}10{\color{Purple})}\^{}2{\color{Purple})} {\color{DarkGreen}/} {\color{Purple}(}$\zeta${\color{Purple}(}5{\color{Purple})}\^{}2 $\zeta${\color{Purple}(}18{\color{Purple})}{\color{Purple})}}}\\
\hspace*{1cm}{\small\texttt{{\color{Blue}[C-25]} L{\color{Purple}(}$\lambda$ $\tau$ $\sigma$', 6{\color{Purple})} {\color{FireBrick}=} {\color{Purple}(}$\zeta${\color{Purple}(}5{\color{Purple})}\^{}2 $\zeta${\color{Purple}(}11{\color{Purple})} $\zeta${\color{Purple}(}12{\color{Purple})}\^{}2{\color{Purple})} {\color{DarkGreen}/} {\color{Purple}(}$\zeta${\color{Purple}(}6{\color{Purple})}\^{}2 $\zeta${\color{Purple}(}22{\color{Purple})}{\color{Purple})}}}\\
\hspace*{1cm}{\small\texttt{{\color{Blue}[C-25]} L{\color{Purple}(}$\lambda$ $\tau$ $\sigma$', 7{\color{Purple})} {\color{FireBrick}=} {\color{Purple}(}$\zeta${\color{Purple}(}6{\color{Purple})}\^{}2 $\zeta${\color{Purple}(}13{\color{Purple})} $\zeta${\color{Purple}(}14{\color{Purple})}\^{}2{\color{Purple})} {\color{DarkGreen}/} {\color{Purple}(}$\zeta${\color{Purple}(}7{\color{Purple})}\^{}2 $\zeta${\color{Purple}(}26{\color{Purple})}{\color{Purple})}}}\\
\hspace*{1cm}{\small\texttt{{\color{FireBrick}[!!!!]} L{\color{Purple}(}$\lambda$ $\tau$, 4{\color{Purple})} {\color{FireBrick}=} $\zeta${\color{Purple}(}8{\color{Purple})}\^{}2 {\color{DarkGreen}/} $\zeta${\color{Purple}(}4{\color{Purple})}\^{}2}}\\
\hspace*{1cm}{\small\texttt{{\color{FireBrick}[!!!!]} L{\color{Purple}(}$\lambda$ $\tau$, 5{\color{Purple})} {\color{FireBrick}=} $\zeta${\color{Purple}(}10{\color{Purple})}\^{}2 {\color{DarkGreen}/} $\zeta${\color{Purple}(}5{\color{Purple})}\^{}2}}\\
\hspace*{1cm}{\small\texttt{{\color{FireBrick}[!!!!]} L{\color{Purple}(}$\lambda$ $\tau$, 6{\color{Purple})} {\color{FireBrick}=} $\zeta${\color{Purple}(}12{\color{Purple})}\^{}2 {\color{DarkGreen}/} $\zeta${\color{Purple}(}6{\color{Purple})}\^{}2}}

Here, the remaining unknown finds are instances of C-22.

\end{document}